\documentclass[12pt]{article}

\usepackage{amsfonts}
\usepackage{amsmath}
\usepackage{esint}
\usepackage[margin=1in]{geometry}
\usepackage{graphicx}
\usepackage{epstopdf}
\usepackage{amssymb}
\usepackage{cite}
\usepackage{multicol}
\usepackage{caption}
\usepackage{color}
\usepackage[section]{placeins}
\usepackage[linesnumbered,ruled,vlined]{algorithm2e}
\usepackage[title]{appendix}
\usepackage{hyperref}

\newcommand{\ds} {\displaystyle}
\newcommand{\Frac}[2]{\ds \frac{#1}{#2}}
\newcommand{\R}{{\mathbb R}}

\SetKwInput{KwInput}{Input}      
\SetKwInput{KwOutput}{Output}

\title {Fourier Series-Based Approximation of Time-Varying Parameters in Ordinary Differential Equations}
\author{Anna Fitzpatrick, Molly Folino, Andrea Arnold*}  
\date{}

\begin{document}
\maketitle

\small 
\centerline{Department of Mathematical Sciences, Worcester Polytechnic Institute, Worcester, MA, USA}  
\vspace{.2cm}  

\centerline{* Corresponding author: anarnold@wpi.edu}  

\normalsize

\bigskip

\begin{abstract}
Many real-world systems modeled using differential equations involve unknown or uncertain parameters. 
Standard approaches to address parameter estimation inverse problems in this setting typically focus on estimating constants; yet some unobservable system parameters may vary with time without known evolution models.
In this work, we propose a novel approximation method inspired by the Fourier series to estimate time-varying parameters in deterministic dynamical systems modeled with ordinary differential equations.  
Using ensemble Kalman filtering in conjunction with Fourier series-based approximation models, we detail two possible implementation schemes for sequentially updating the time-varying parameter estimates given noisy observations of the system states. 
We demonstrate the capabilities of the proposed approach in estimating periodic parameters, both when the period is known and unknown, as well as non-periodic time-varying parameters of different forms with several computed examples using a forced harmonic oscillator.
Results emphasize the importance of the frequencies and number of approximation model terms on the time-varying parameter estimates and corresponding dynamical system predictions. \\

\noindent \textbf{Keywords:} Nonstationary inverse problems; parameter estimation; approximation models; Bayesian inference; Fourier series; ensemble Kalman filter; dynamical systems.
\end{abstract}



\section{Introduction}

Many real-world problems in science and engineering involve unknown model parameters of interest that may vary with time but cannot be directly observed.
For example, forced harmonic oscillators used in modeling gear systems \cite{Sika2008, Shen2006}, RLC circuits \cite{Batouli2015}, and cantilever motion in atomic force microscopy \cite{Preiner2007} may involve time-dependent stiffness, mass, and/or external forcing parameters. 
Examples in biology and medicine include time-dependent transmission parameters in modeling epidemic dynamics \cite{Altizer2006, Zeng2020, Calvetti2020}, 
external stimuli in modeling neuron dynamics \cite{Linaro2018, Shamir2009, Campbell2020}, 
and tissue optical properties in modeling laser-tissue interactions \cite{ArnoldFichera2022, Bashkatov2016}.
The work in this paper aims to address the estimation of such time-varying parameters (TVPs) in deterministic dynamical systems, when the parameters are not observable and there is not a known (or available) model governing their time evolution.

While TVPs may appear in different types of mechanistic models, here we focus on deterministic dynamical systems modeled using differential equations.
In particular, we assume an ordinary differential equation (ODE) model of the form
\begin{equation}\label{Eq:State_Dynamics}
\frac{dx}{dt} = f(t,x,\theta), \quad x(0) = x_0
\end{equation}
where $t\in\R$ denotes time, $x=x(t)\in\R^d$ is the vector of model states, $\theta\in\R^p$ is the vector of unknown model parameters, and $f:\R\times\R^d\times\R^p\rightarrow\R^d$ is a known mapping representing the state dynamics.  
The ODE model in \eqref{Eq:State_Dynamics} may involve unknown constant parameters (including the initial conditions $x_0$), but here we focus our attention primarily on estimating the unknown time-varying parameters of interest.  
Further, while such problems may generally include more than one parameter changing with time, we restrict our examples in this work to estimate a single, univariate TVP, $\theta = \theta(t)\in\R$, for each system.
The inverse problem considered is therefore to estimate $\theta(t)$, along with the system states $x(t)$, at some discrete times $t_j$ given noisy, sequential observations of the system states (which may be fully or partially observed). 
Note that we can extend the inverse problem to include estimation of additional constant parameters, including initial conditions, as needed.

To address the inverse problem at hand, in this work we propose a novel approximation method inspired by the Fourier series, where we represent $\theta(t)$ as a linear combination of a finite number of sine and cosine functions and estimate the unknown coefficients of the approximation model using ensemble Kalman filtering.  
Fourier series-based approaches have been used in previous work for approximating smooth, periodic functions and trajectories in control systems \cite{Labiod2013, Caruso2021} and for linear function approximation in reinforcement learning \cite{Konidaris2011}. 
Fourier series expansion has also been used for estimating smooth, periodic TVPs in the setting of adaptive control \cite{Liuzzo2007, Chen2010, Zhang2012, Chen2021}, where the period of the parameter is assumed to be known.
Our goal in this work is to approximate more general $\theta(t)$, specifically periodic parameters for which the period is not known and parameters that are time-varying but not periodic over the time interval of available system observations.

Use of a Bayesian filtering approach, such as ensemble Kalman filtering, provides a flexible framework for addressing this inverse problem and allows for the unknown approximation model coefficients to be updated sequentially along with the system states as new data arrive, without relying on the full time series of data in advance.
Ensemble-based methods also provide a natural measure of uncertainty in the resulting parameter estimates, generally taken as $\pm 2$ standard deviations around the sample mean at each time.  
In the setting of TVP estimation, previous studies have used Bayesian filtering approaches to estimate the constant coefficients of piecewise functional representations of periodic parameters \cite{Arnold2018, Arnold2020}, constant parameters for linear TVP models \cite{Feng2017}, and in conjunction with parameter tracking schemes to estimate more general time-varying forms \cite{Voss2004, Bian2011, Arnold2019, Campbell2020, Calvetti2021, Arnold2023}.

In a set of numerical examples utilizing a forced mass-spring system, we first establish the ability of the proposed Fourier series-based approximation approach to estimate a periodic forcing parameter with known period, then extend the estimation to include the period of the TVP as an additional unknown.
We further demonstrate the capability of the proposed method in estimating TVPs that are not periodic (and potentially not continuous) over the time interval of observed data and show that the resulting parameter approximation models can be used to make reasonably accurate predictions of the system dynamics for different initial conditions.

The remainder of the paper is organized as follows:
Section~\ref{Sec:EnKF} briefly reviews ensemble Kalman filtering for constant parameter estimation.
Section~\ref{Sec:Fourier} details the proposed Fourier series-based TVP approximation models, providing two implementation approaches that can be used together with ensemble Kalman filtering to estimate the unknown model coefficients.
Section~\ref{Sec:Results} provides the main results of the numerical experiments (with some additional results provided in Appendix~\ref{Appendix:Augmentation_Approach}), and Section~\ref{Sec:Conclusions} gives conclusions and future work.



\section{Review: Ensemble Kalman Filtering for Sequential Estimation of Constant Parameters}
\label{Sec:EnKF}

\begin{figure}[t!]
\centerline{\includegraphics[width = 0.65\textwidth]{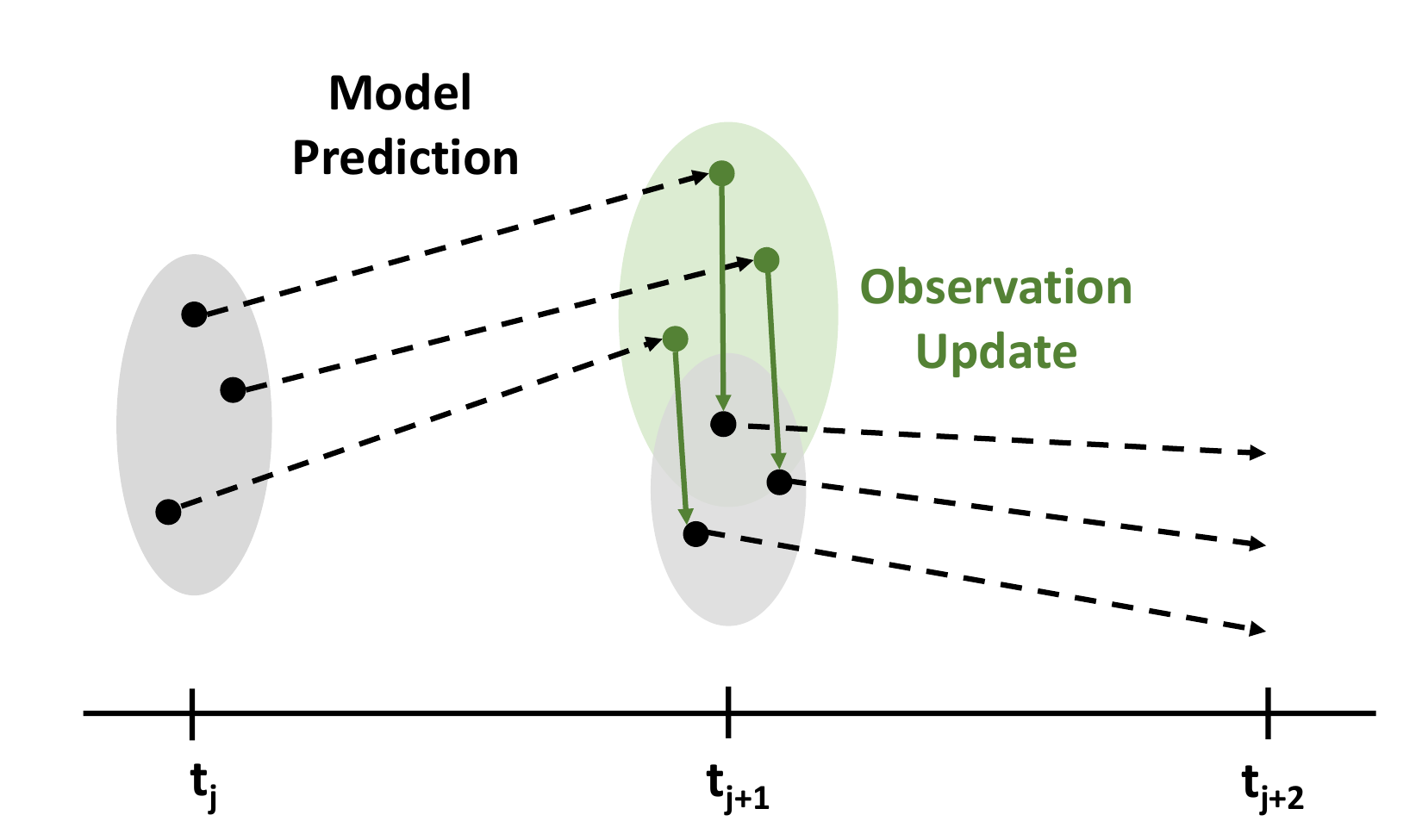} }
\caption{Illustration of the EnKF two-step updating scheme.  In the model prediction step, the ensemble members (represented as black circles, with sample variance in gray) at time $j$ are propagated forward (shown with black dashed arrows) to time $j+1$ by solving the model in \eqref{Eq:Aug_Dynamics}.  In the observation update, the ensemble predictions (green circles) are corrected (green arrows) using the observed system data at time $j+1$.  The process continues using the corrected sample at time $j+1$ (black circles) until all available data are assimilated. }
\label{Fig:EnKF}
\end{figure}

In this section, we briefly review the main ideas behind ensemble Kalman filtering for combined state and constant parameter estimation.
While the method was originally established for tracking unobserved model states \cite{Evensen1994, Burgers1998}, the augmented approach outlined below has been successfully used to estimate constant parameters in a variety of models, including systems of ODEs \cite{Arnold2014}.
For a recent review of ensemble Kalman filtering, we refer interested readers to \cite{Katzfuss2016}.

The Ensemble Kalman Filter (EnKF) is a sequential Bayesian approach that employs ensemble statistics within the framework of the classic Kalman filter to track unknown model variables given observed time series data.  
The incorporation of a statistical sample, which represents an underlying probability distribution of the unknowns conditioned on the available data, accommodates the use of nonlinear and possibly non-Gaussian models.  
The EnKF can estimate unobserved (or unobservable) model states along with unknown constant parameters through use of an augmented dynamical system
\begin{equation}\label{Eq:Aug_Evolution}
\Frac{dz}{dt} = \left[ \begin{array}{c} \Frac{dx}{dt} \\[.4cm] \Frac{d\theta}{dt} \end{array} \right] , \quad z(0) = \left[ \begin{array}{c} x(0) \\[.1cm] \theta(0) \end{array} \right] \in \R^{d+p}
\end{equation}
where $dx/dt$ describes the system dynamics, as given in \eqref{Eq:State_Dynamics}, while $d\theta/dt$ represents the dynamics of the parameters \cite{Evensen2009, Arnold2014}.  
In particular, $d\theta/dt = 0$ when the parameters are constant (i.e., time-invariant).

Given a discrete sample of the model states and parameters at time $j$,
\begin{equation}
\mathcal{S}_j = \big\{ (x_j^{(1)}, \theta_j^{(1)}), \dots, (x_j^{(N)}, \theta_j^{(N)}) \big\} 
\end{equation}
the EnKF works as a two-step updating scheme: First, each pair of states and parameters is predicted at time $j+1$ using the evolution model in \eqref{Eq:Aug_Evolution}; then, the augmented vectors are corrected using the Kalman filter observation updating equation, which incorporates the observed data at time $j+1$.  
This process is illustrated in Figure~\ref{Fig:EnKF}, and the steps of the algorithm for combined state and constant parameter estimation are outlined in Algorithm~\ref{Alg:EnKF}.

\begin{algorithm}[t!]
    \DontPrintSemicolon
    \caption{EnKF for Combined State and Constant Parameter Estimation \label{Alg:EnKF}}
    \KwInput{Initial sample $\mathcal{S}_0$ drawn from prior distribution $\pi(x_0,\theta_0)$}
    \KwOutput{Posterior sample $\mathcal{S}_T$ and corresponding ensemble statistics} 
    
    Initialize time index $j = 0$
    
    \While{$j < T$}
    {
        
    \tcc{Prediction Step}
    \For{n = 1, \dots, N}    
    { 
    $x_{j+1}^{(n)} = F(x_{j}^{(n)},\theta_{j}^{(n)}) + v_{j+1}^{(n)}, \quad v_{j+1}^{(n)}\sim\mathcal{N}(0,\mathsf{C})$  
    
    $z_{j+1}^{(n)} = \big[ x_{j+1}^{(n)}; \theta_j^{(n)} \big]$
    } 
        
    \tcc{Observation Update}
    \For{n = 1, \dots, N}    
    { 
    $y_{j+1}^{(n)} = y_{j+1} + w_{j+1}^{(n)}, \quad w_{j+1}^{(n)}\sim\mathcal{N}(0,\mathsf{D})$  
    
    $z_{j+1}^{(n)} = z_{j+1}^{(n)} + \mathsf{K}_{j+1} \big( y_{j+1}^{(n)} - G(z_{j+1}^{(n)}) \big)$  
    }
    
    \tcc{Compute Posterior Ensemble Statistics}
    
    $\bar{z}_{j+1} = \frac{1}{N} \sum_{n=1}^{N} z_{j+1}^{(n)}$
       
    $\mathsf{\Gamma}_{j+1} = \frac{1}{N-1} \sum_{n=1}^{N} (z_{j+1}^{(n)} - \bar{z}_{j+1})(z_{j+1}^{(n)} - \bar{z}_{j+1})^\mathsf{T}$
    
    \tcc{Update Time Index}
    $j = j+1$
    }
 \end{algorithm}

When working with ODE models, the operator $F$ in the prediction step of Algorithm~\ref{Alg:EnKF} (line 4) denotes the numerical solution to the differential equations model in \eqref{Eq:State_Dynamics} at time $j+1$; the parameter values are not updated during this step.  
In the observation update, $y_{j+1}$ in line 7 denotes the vector of observed model states, with dimension $m\leq d$, which is perturbed to help avoid too low a covariance in the resulting sample \cite{Burgers1998}.  
The updating equation in line 8 corrects the joint predicted sample for each $n$ using the Kalman gain matrix $\mathsf{K}_{j+1}$ and the difference between the perturbed observation $y_{j+1}^{(n)}$ and predicted observation $G(z_{j+1}^{(n)})$, where $G$ is the observation model.  
For linear observations, $G(z_{j+1}^{(n)}) = \mathsf{P} z_{j+1}^{(n)}$, where $\mathsf{P}\in\R^{m\times(d+p)}$ is a projection matrix whose entries corresponding to observed states are 1 and entries corresponding to unobserved states and parameters are 0.  
The posterior ensemble statistics computed in lines 9 and 10 give the mean and covariance, respectively, of the resulting sample at time $j+1$.  
The posterior mean for each parameter is taken as its estimate, with uncertainty commonly represented using $\pm2$ standard deviations around the mean.
The process repeats sequentially until all available data in the time series are assimilated.  
In this procedure, the parameter values are artificially evolved with the aim of converging to a constant.



\section{Fourier Series-Based Approximation Models for Time-Varying Parameters}
\label{Sec:Fourier}

While the EnKF algorithm reviewed in Section~\ref{Sec:EnKF} is formulated for estimating constant parameters, our goal in this work is to estimate time-varying system parameters for which no evolution model is known (or available); i.e., we do not have a known form of $d\theta/dt$ for $\theta = \theta(t)$.  
To address this problem, we propose a novel approximation method inspired by the Fourier series, where $\theta(t)$ is represented as a linear combination of sine and cosine functions with different frequencies.  
We describe this approach below, detailing formulations for estimating periodic parameters when the period of $\theta(t)$ is both known and unknown, as well as for approximating more general TVPs that are not periodic over the time interval of observed data.
We provide two possible implementation strategies using the EnKF for coefficient estimation.

\subsection{Fourier Series and Approximation Model Formulation}
\label{Sec:Approx_Model}

\begin{figure}[t!]
\centerline{\includegraphics[width = \textwidth]{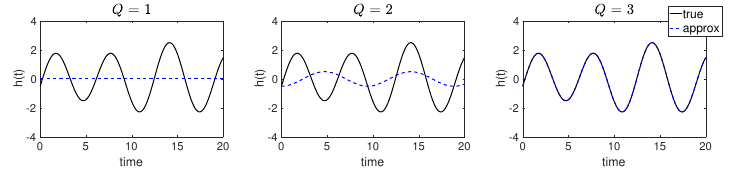} }
\caption{Fourier series approximations of the sinusoidal function $h(t) = 2\sin(t) - 0.5\cos(2t/3)$ using different values of $Q$, which dictates the number of terms used in the approximation.
In each plot, the Fourier series approximation is shown in dashed blue, while the true function is shown in solid black. 
The integrals in \eqref{Eq:Fourier_acoeff} and \eqref{Eq:Fourier_bcoeff} were computed numerically using adaptive quadrature via the \texttt{integral} function in MATLAB \cite{Shampine2008}.
}
\label{Fig:FourierEx}
\end{figure}

Recall that the $Q$th order Fourier series of a univariate function $h(t)$ is given by
\begin{equation}\label{Eq:Fourier_series}
h(t) = \frac{a_0}{2} + \sum_{q=1}^{Q} \bigg( a_q \cos\Big(\frac{2\pi qt}{P}\Big) + b_q\sin\Big(\frac{2\pi qt}{P}\Big) \bigg)
\end{equation}
where $a_q$ and $b_q$ are the expansion coefficients and $P$ is the period.  If $h(t)$ is known, the coefficients $a_q$ and $b_q$ can be computed explicitly using the formulas
\begin{equation}\label{Eq:Fourier_acoeff} 
a_q = \Frac{2}{P} \ds\int_0^P h(t) \cos\Big(\frac{2\pi qt}{P}\Big) dt 
\end{equation}
and
\begin{equation}\label{Eq:Fourier_bcoeff}
b_q = \Frac{2}{P} \ds\int_0^P h(t) \sin\Big(\frac{2\pi qt}{P}\Big) dt 
\end{equation}
respectively.
For example, Figure~\ref{Fig:FourierEx} shows the Fourier series approximations of the periodic function $h(t) = 2\sin(t) - 0.5\cos(2t/3)$ with known period $P=6\pi$ for different choices of $Q$.
However, when using Fourier series for function approximation, the function $h(t)$ is generally unknown and the coefficients $a_q$ and $b_q$ must be estimated \cite{Konidaris2011, Liuzzo2007, Zhang2012}.

Inspired by use of the Fourier series for approximating unknown functions, we propose to represent the unknown, univariate time-varying parameters $\theta(t)$ in this work as linear combinations of sine and cosine pairs, such that the estimate of each $\theta(t)$ is determined by the approximation model
\begin{equation} \label{Eq:TV_form}
\theta_M(t) = \alpha_0 + \sum_{i=1}^M \Big( \alpha_i \sin(\omega_i t) + \beta_i \cos(\omega_i t) \Big)
\end{equation}
with coefficients $\alpha_i$, $\beta_i$ and fixed values of $\omega_i$ for each $i=1,\dots,M$.  
Moving forward, we denote the $2M+1$ coefficients by $c_k$, $k=0,\dots,2M$, such that the model in \eqref{Eq:TV_form} becomes 
\begin{equation} \label{Eq:TV_form2}
\theta_M(t) = c_0 + c_1 \sin(\omega_1 t) + c_2 \cos(\omega_1 t) + \cdots + c_{2M-1} \sin(\omega_M t) + c_{2M} \cos(\omega_M t)
\end{equation}
for some fixed $M$.
Our goal in approximating $\theta(t)$ is therefore is to estimate the $2M+1$ coefficients $c_k$ that provide the best fit between the model in \eqref{Eq:TV_form2} and the true time-varying parameter given the observed system data.  
This essentially transforms the inverse problem at hand into a constant parameter estimation problem, and we can utilize the EnKF algorithm described in Section~\ref{Sec:EnKF} to estimate these coefficients.

To apply the TVP approximation model in \eqref{Eq:TV_form2}, one must specify the value of $M$, which sets the number of terms in the approximation, and the values of $\omega_i$, $i=1,\dots,M$, which act as the angular frequencies of the sinusoidal functions.
When $\theta(t)$ is periodic, we consider two different approaches in assigning the $\omega_i$ values: 
If we know the period of the underlying $\theta(t)$ in advance of the estimation process, we explicitly define $\omega_i$ using the formula
\begin{equation}\label{Eq:Known_Omega}
\omega_i = \frac{2\pi i}{P}, \quad i = 1,\dots, M
\end{equation}
where $P$ is the known period of $\theta(t)$; this form of $\omega_i$ follows from the terms in the Fourier series in \eqref{Eq:Fourier_series} for approximating a periodic function with known period.
If the period of $\theta(t)$ is unknown or uncertain, we set $\omega_i$ as in \eqref{Eq:Known_Omega} and treat the period $P$ as an additional unknown constant parameter to be jointly estimated along with the coefficients of the TVP approximation model.
In more general cases, i.e., when $\theta(t)$ is not known to be periodic over the time interval of available data, we instead choose a fixed increment $\omega$ and define $\omega_i = \omega i$ for each $i = 1,\dots, M$, systematically incorporating sine and cosine pairs with different frequencies into the approximation.
Following this approach, appropriate choice of the increment $\omega$ becomes an important factor in the estimation process.

\subsection{Alternative Implementation: Derivative-Based Augmentation of System States}
\label{Sec:Alternative}

In the previous section, we propose an approximation model for $\theta(t)$ that we can use directly within the ODE model in \eqref{Eq:State_Dynamics} in place of the TVP.
As an alternative means of implementation, much in the spirit of state augmentation, we can use the approximation model in \eqref{Eq:TV_form2} to prescribe a model for $d\theta/dt$ and define a coupled ODE system of the form
\begin{eqnarray}
\frac{dx}{dt} &=& f(t,x,\theta), \quad x(0) = x_0 \\[0.2cm]
\frac{d\theta}{dt} &=& \theta'_M(t), \quad \theta(0) = \theta_0
\end{eqnarray}
which we can write equivalently as the augmented system
\begin{equation}\label{Eq:Aug_Dynamics}
\Frac{dz}{dt} = \left[ \begin{array}{c} \Frac{dx}{dt} \\[.4cm] \Frac{d\theta}{dt} \end{array} \right] = \left[ \begin{array}{c} f(t,x,\theta) \\[.2cm] \theta'_M(t) \end{array} \right] , \quad z(0) = z_0 = \left[ \begin{array}{c} x_0 \\[.2cm] \theta_0 \end{array} \right] 
\end{equation}
where $z = z(t) \in \R^{d+1}$.
In this representation, $\theta = \theta(t)$ is treated as an unobserved state of the system in \eqref{Eq:Aug_Dynamics}.
The equation for $d\theta/dt$ in \eqref{Eq:Aug_Dynamics} follows from taking the derivative of the approximation model $\theta_M(t)$ in \eqref{Eq:TV_form2}, which gives
\begin{equation} \label{Eq:TV_deriv}
\theta'_M(t) = c_1 \omega_1 \cos(\omega_1 t) - c_2 \omega_1 \sin(\omega_1 t) + \cdots + c_{2M-1} \omega_M \cos(\omega_M t) - c_{2M} \omega_M \sin(\omega_M t)
\end{equation}
with $2M$ unknown coefficients, $c_1, \dots, c_{2M}$.
Note that the additive constant $c_0$ no longer explicitly appears in the system equations; instead, we estimate the initial condition $\theta_0$ as an additional constant parameter playing a similar role.



\section{Numerical Results}
\label{Sec:Results}

In this section, we detail the results of several numerical experiments demonstrating the effectiveness of the proposed methodology under different scenarios for $\theta(t)$; more specifically, we consider examples where $\theta(t)$ is a periodic TVP, in cases assuming both a known period and an unknown period, and where $\theta(t)$ is a non-periodic TVP over the time interval of available system data.  
Results were obtained using MATLAB$^{\scriptsize\textregistered}$ (The MathWorks, Inc., Natick, MA) programming language.  


As a test system in the computed examples that follow,
we consider a forced harmonic oscillator, classically modeled using the second-order ODE 
\begin{equation}\label{eq:SpringMass}
mp''+bp'+kp=\theta(t)
\end{equation}
where $p = p(t)$ commonly denotes the position (or displacement) of a mass at time $t$, $m>0$ is the constant mass, $b>0$ is the damping coefficient, and $k>0$ is the spring constant; see, e.g., \cite{Nagle2011, Boyce2001}.  
Here $\theta(t)$ represents external forcing applied to the system, which we treat as our time-varying parameter of interest.  
Letting $v(t) = p'(t)$ denote the velocity of the mass, we can rewrite \eqref{eq:SpringMass} as a first-order ODE system of the form  
\begin{eqnarray}
\frac{dp}{dt} &=& v\\
\frac{dv}{dt} &=& \frac{1}{m}\Big(-kp - bv + \theta(t) \Big)
\end{eqnarray}
or, equivalently, as
\begin{equation}\label{Eq:MS_system}
\frac{dx}{dt} = \left[ \begin{array}{cc} \ 0 & \ 1 \\[.2cm] -\Frac{k}{m} & -\Frac{b}{m} \end{array} \right] x + \left[ \begin{array}{c} 0 \\[.2cm] \Frac{\theta(t)}{m} \end{array} \right] 
\end{equation}
where $x(t) = [p(t); v(t)] \in\R^2$ is the vector of model states at time $t$.  
Assuming that the constants $m$, $k$, and $b$ are known, our goal is to estimate $\theta(t)$ using the Fourier series-based approximation methods described in Section~\ref{Sec:Fourier}.

To test the effectiveness of the proposed estimation techniques, we generate data from the mass-spring system in \eqref{Eq:MS_system} using the initial condition $x(0) = [2; 0]$, fixed constants $m=10$, $k=5$, and $b=3$, and different forms of the time-varying forcing parameter $\theta(t)$, as detailed for each experiment below.
We initialize the EnKF with a sample size of $N = 100$ and draw the prior ensemble of state values from a multivariate Gaussian distribution with mean $[1; 1]\in\R^2$ and covariance matrix $(0.5)^2 \mathsf{I}_2$, where $\mathsf{I}_2$ denotes the $2\times2$ identity matrix.  
We draw the prior ensemble of values for each of the unknown coefficients $c_k$ uniformly over the interval $[-2, 10]$.  
Further, we use MATLAB's \texttt{ode15s} to solve the ODE system in \eqref{Eq:MS_system} at each time step of the filter and prescribe $\mathsf{C} = (0.02)^2 \mathsf{I}_2$ as the model innovation covariance matrix and $\mathsf{D} = (0.08)^2 \mathsf{I}_2$ as the observation covariance matrix, assuming observations of both position and velocity.


\subsection{Example: Periodic Time-Varying Parameter}

\begin{figure}[t!] 
\centerline{\includegraphics[width = \textwidth]{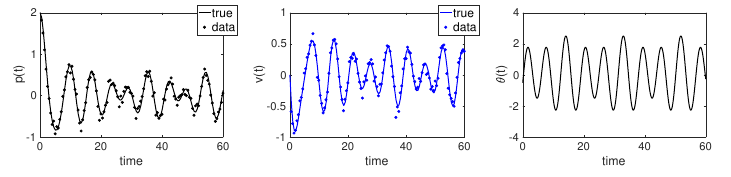}}
\caption{Position and velocity data generated from the mass-spring system in \eqref{Eq:MS_system} with sinusoidal forcing parameter $\theta(t) = 2\sin(t) - 0.5\cos(2t/3)$. }
\label{Fig:MS_Data}
\end{figure}

In this example, we consider a sinusoidal forcing parameter of the form $\theta(t) = 2\sin(t) - 0.5\cos(2t/3)$ with period $P=6\pi$ as the underlying truth.  
Utilizing MATLAB's \texttt{ode45} to solve the ODE system in \eqref{Eq:MS_system}, we record observations every 0.5 time units over the interval [0,60], spanning just over three periods of $\theta(t)$, and corrupt the observations using Gaussian noise with zero mean and standard deviation taken to be 20\% of the standard deviation of the true system states.  
Figure~\ref{Fig:MS_Data} shows the simulated data.  
In the experiments that follow, we consider two cases for the estimation procedure: (i) when the period of $\theta(t)$ is known, and (ii) when the period of $\theta(t)$ is unknown.
Results focus on use of the TVP approximation model approach described in Section~\ref{Sec:Approx_Model}.
Appendix~\ref{Appendix:Augmentation_Approach} shows an example of the corresponding numerical results when using the derivative-based augmentation approach described in Section~\ref{Sec:Alternative} for the known period case.

\subsubsection{Estimation with Known Period} 
\label{Sec:KnownPeriod_Ex}

\begin{figure}[t!] 
\centerline{\includegraphics[width = \textwidth]{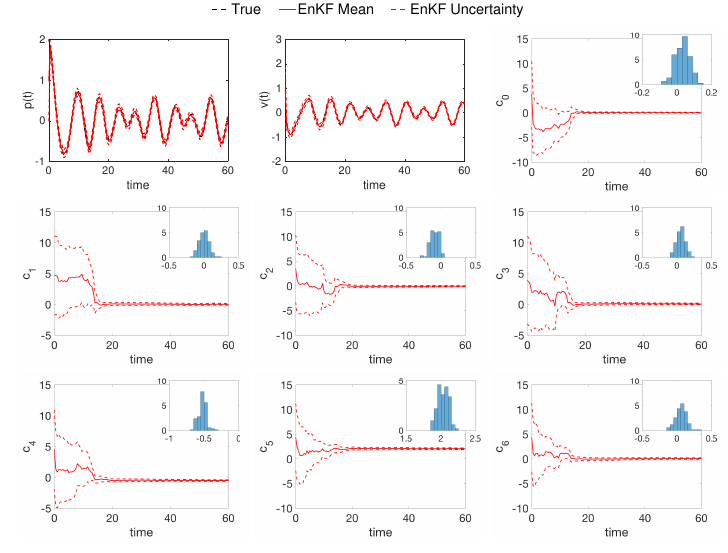}}
\caption{EnKF time series estimates of position, velocity, and the seven unknown coefficients in the TVP approximation model \eqref{Eq:TV_form2} of the sinusoidal forcing parameter $\theta(t) = 2\sin(t) - 0.5\cos(2t/3)$ in \eqref{Eq:MS_system} when $M=3$, using the data in Figure~\ref{Fig:MS_Data} and assuming a known period for $\theta(t)$.
On each plot, the solid red line denotes the EnKF sample mean and the dashed red lines show $\pm2$ standard deviations around the mean.
The plot for each coefficient, $c_0, \dots, c_6$, also shows the resulting histogram of the posterior sample at time $t=60$.
}
\label{Fig:KnownPeriod_Results_M3}
\end{figure}

\begin{figure}[t!] 
\centerline{\includegraphics[width = \textwidth]{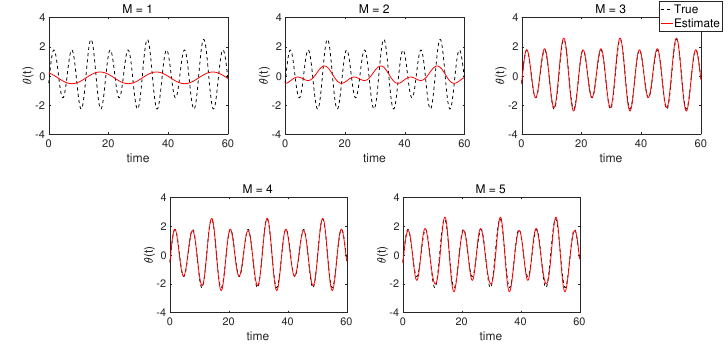}}
\caption{Resulting Fourier series-based approximations of the sinusoidal forcing parameter $\theta(t) = 2\sin(t) - 0.5\cos(2t/3)$ in \eqref{Eq:MS_system} for different values of $M$, using the data in Figure~\ref{Fig:MS_Data} and assuming a known period for $\theta(t)$.
On each plot, the solid red line denotes the TVP approximation $\bar{\theta}_{M}(t)$ computed using the posterior mean coefficient values in Table~\ref{Tab:KnownPeriod_Results} and the dashed black line shows the true underlying $\theta(t)$.
}
\label{Fig:KnownPeriod_Results_TVP}
\end{figure}

Assuming a known period of $P=6\pi$ for $\theta(t)$, we apply the formula in \eqref{Eq:Known_Omega} to set $\omega_i = {i}/{3}$, $i = 1,\dots, M$,
and employ the TVP approximation model in \eqref{Eq:TV_form2} with increasing values of $M$, 
ranging from $M=1$ (which involves three coefficients) through $M=5$ (11 coefficients).  
In each case, the EnKF tracks the model states and estimates the unknown coefficients comprising the Fourier series-based model approximation of $\theta(t)$.
As an example, Figure~\ref{Fig:KnownPeriod_Results_M3} shows the resulting EnKF time series estimates of position, velocity, and the seven TVP model coefficients when $M=3$.
We then use the posterior sample mean of each coefficient, $\bar{c}_k$, $k=0,\dots,2M$, to construct an approximation of $\theta(t)$, such that 
\begin{equation} \label{Eq:TVP_EnKF}
\bar{\theta}_{M}(t) = \bar{c}_0 + \bar{c}_1 \sin(\omega_1 t) + \bar{c}_2 \cos(\omega_1 t) + \cdots + \bar{c}_{2M-1} \sin(\omega_M t) + \bar{c}_{2M} \cos(\omega_M t)
\end{equation}
for each $M$.
Table~\ref{Tab:KnownPeriod_Results} lists the posterior sample mean of each coefficient for each $M$ value, and
Figure~\ref{Fig:KnownPeriod_Results_TVP} shows the resulting $\bar{\theta}_{M}(t)$ approximations in each case compared with the true $\theta(t)$. 
To further compare the approximations with the true $\theta(t)$, we use a scaled version of the root mean square error (RMSE), where
\begin{equation}
\text{Scaled RMSE} = \frac{\text{RMSE}}{\sigma_\theta}, \qquad \text{RMSE} = \sqrt{\frac{1}{T}\sum_{j=1}^T \big( \theta(t_j) - \bar{\theta}_{M}(t_j) \big)^2}
\end{equation}
and $\sigma_\theta$ is the standard deviation of $\theta(t)$.
Note that the plots in Figure~\ref{Fig:KnownPeriod_Results_TVP} and corresponding scaled RMSE values in Table~\ref{Tab:KnownPeriod_Results} were computed at $t_j$ values taken every 0.1 time units over [0,60], a finer time discretization than used during the filtering process.

As illustrated in Figure~\ref{Fig:KnownPeriod_Results_M3}, the filter well tracks the model states for both $p(t)$ and $v(t)$, and the coefficient estimates converge to constant values after assimilating approximately one period of data (here, one period occurs at time $6\pi \approx 18.85$), with the $\pm2$ standard deviation curves representing uncertainty around the mean estimate shrinking significantly.  
In comparing the results for different values of $M$, we note that the lowest RMSE occurs when $M=3$ and is similar when $M=4$, with a small increase in error when $M = 5$; however, all three of these $M$ values result in reasonably close model approximations to the true $\theta(t)$, as shown in Figure~\ref{Fig:KnownPeriod_Results_TVP}.
Further, the EnKF posterior sample means for the TVP approximation model coefficients when $M=3$, 4, and 5 are quite similar to the Fourier series coefficients obtained when using the formulas in \eqref{Eq:Fourier_acoeff} and \eqref{Eq:Fourier_bcoeff} to approximate the function $h(t) = 2\sin(t) - 0.5\cos(2t/3)$ when $Q=3$, 4, and 5, assuming that both $h(t)$ and $P=6\pi$ are known.

\begin{figure}[t!] 
\centerline{\includegraphics[width = \textwidth]{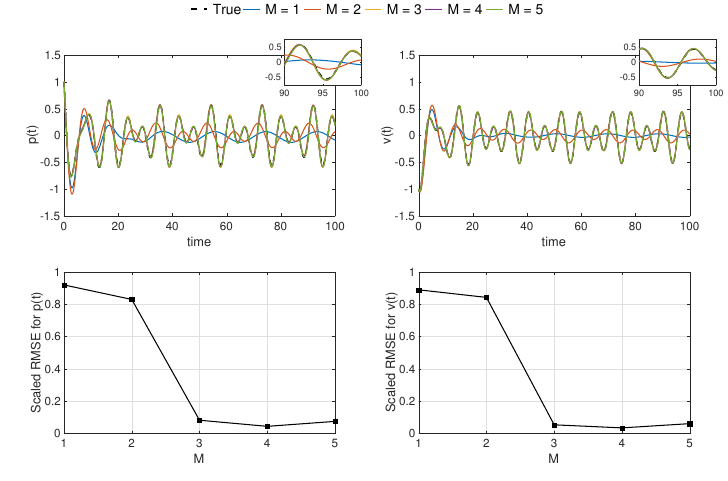}}
\caption{Top row: Predictions of position (left) and velocity (right) using initial condition $p(0) = 1$, $v(0) = -1$ and the TVP approximations for $\theta(t)$ computed with the estimated coefficients in Table~\ref{Tab:KnownPeriod_Results} for different $M$.  Bottom row: Scaled RMSE comparing the predictions of position (left) and velocity (right) with the true solutions for each $M$. 
}
\label{Fig:KnownPeriod_Pred}
\end{figure}

Figure~\ref{Fig:KnownPeriod_Pred} shows the mass-spring system model predictions using $\bar{\theta}_{M}(t)$ as the forcing parameter in \eqref{Eq:MS_system} for each $M$ with the initial condition $x(0) = [1; -1]$, a different initial condition than used in generating the simulated data, along with the corresponding scaled RMSE values for each model state.
As might be expected given the TVP approximation results, the model predictions when using $M=1$ and 2 are less accurate than when using $M=3$, 4, and 5, which all provide similarly accurate predictions of both position and velocity.
While not shown, similar results hold for other initial conditions in this range.
For predictions using initial conditions larger in magnitude with this system, e.g., $x(0) = [100; -100]$, a similar pattern holds where the approximation models with $M=3$, 4, and 5 provide the best results (i.e., lowest RMSE values), but the overall error is lower for predictions using any of the $M$ values considered.

\begin{table}[t!]
\caption{EnKF posterior sample means of the TVP model coefficients for different values of $M$ when approximating the sinusoidal forcing parameter $\theta(t) = 2\sin(t) - 0.5\cos(2t/3)$ in \eqref{Eq:MS_system}, using the data in Figure~\ref{Fig:MS_Data} and assuming a known period for $\theta(t)$. Numbers are rounded to four decimal places.}
\label{Tab:KnownPeriod_Results}
\centering
\def\arraystretch{1.25}
\begin{tabular}{cccccc}
Coefficient & $M=1$ & $M=2$ & $M=3$ & $M=4$ & $M=5$ \\[0.2mm]
\hline 
$c_0$ & -0.1328  & -0.0306 & \ 0.0303 & -0.0083 & -0.0552 \\[0.2mm]
$c_1$ & -0.2106  & -0.3373 & -0.0013 & -0.0000 & \ 0.0255 \\[0.2mm]
$c_2$ & \ 0.3387 & -0.2012 & -0.0927 & -0.0624 & -0.0784 \\[0.2mm]
$c_3$ & -- & \ 0.1689 & \ 0.0545 & \ 0.0195 & -0.0035 \\[0.2mm]
$c_4$ & -- & -0.2770 & -0.5236 & -0.5032 & -0.4689 \\[0.2mm]
$c_5$ & -- & -- & \ 2.0299 & \ 2.0047 & \ 2.0538 \\[0.2mm]
$c_6$ & -- & -- & \ 0.0594 & -0.0032 & -0.0533 \\[0.2mm]
$c_7$ & -- & -- & -- & \ 0.1099 & \ 0.1303 \\[0.2mm]
$c_8$ & -- & -- & -- & \ 0.0431 & \ 0.0053 \\[0.2mm]
$c_9$ & -- & -- & -- & -- & -0.1584 \\[0.2mm]
$c_{10}$ & -- & -- & -- & -- & \ 0.1012 \\[0.2mm]
\hline
Scaled RMSE & \ 1.0211  & \ 1.0087  & \ 0.0645  & \ 0.0657  & \ 0.1282
\end{tabular}
\end{table}

\subsubsection{Estimation with Unknown Period} 

\begin{figure}[t!] 
\centerline{\includegraphics[width = \textwidth]{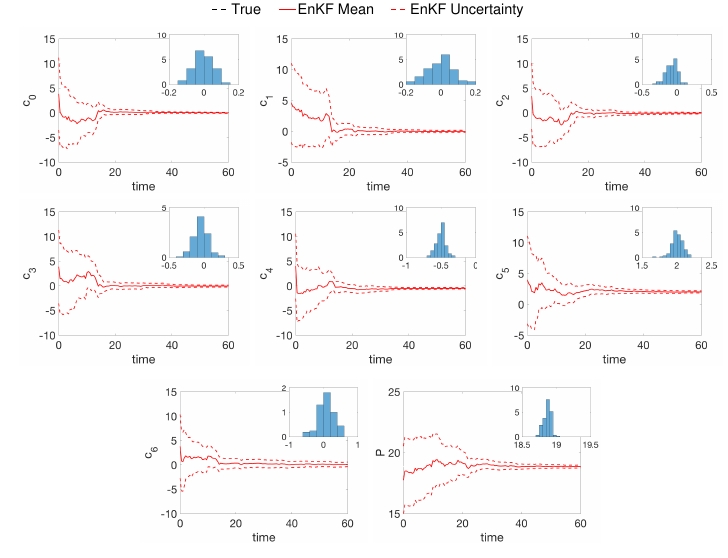}}
\caption{EnKF time series estimates of the seven unknown coefficients and period $P$ in the TVP approximation model \eqref{Eq:TV_form2} of the sinusoidal forcing parameter $\theta(t) = 2\sin(t) - 0.5\cos(2t/3)$ in \eqref{Eq:MS_system} when $M=3$, using the data in Figure~\ref{Fig:MS_Data} and treating the period of $\theta(t)$ as an additional unknown.
On each plot, the solid red line denotes the EnKF sample mean and the dashed red lines show $\pm2$ standard deviations around the mean.
Each plot also shows the resulting histogram of the posterior sample at time $t=60$.
}
\label{Fig:UnknownPeriod_Results_M3}
\end{figure}

\begin{figure}[t!] 
\centerline{\includegraphics[width = \textwidth]{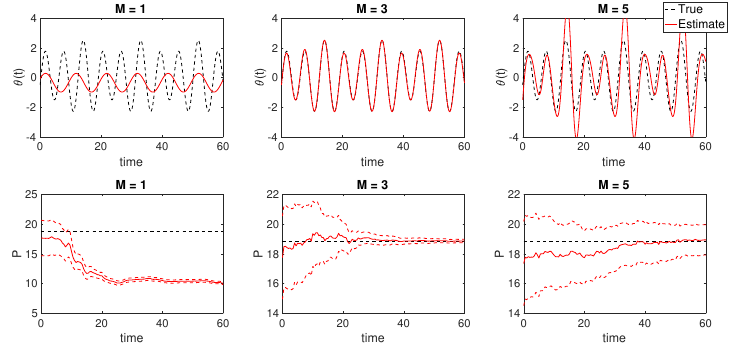}}
\caption{Resulting Fourier series-based approximations of the sinusoidal forcing parameter $\theta(t) = 2\sin(t) - 0.5\cos(2t/3)$ in \eqref{Eq:MS_system} and corresponding time series estimates of the period $P$ for different values of $M$, using the data in Figure~\ref{Fig:MS_Data} and treating the period of $\theta(t)$ as an additional unknown.
Top row: On each plot, the solid red line denotes the TVP approximation $\bar{\theta}_{M}(t)$ computed using the posterior mean coefficient values in Table~\ref{Tab:UnknownPeriod_Results} and the dashed black line shows the true underlying $\theta(t)$.
Bottom row: On each plot, the solid red line denotes the EnKF sample mean, the dashed red lines show $\pm2$ standard deviations around the mean, and the dashed black line shows the true underlying $P$.
}
\label{Fig:UnknownPeriod_Results_TVP}
\end{figure}

Given the same simulated data, we now assume that the period $P$ of $\theta(t)$ is also unknown and estimate $P$ along with the unknown approximation model coefficients for different choices of $M$.
We draw an initial sample of $P$ values from a uniform distribution over [15, 20], supposing that we have some reasonable prior information on a range of likely values (the true period being $6\pi\approx 18.8496$ in this case).
As an example for comparison, Figure~\ref{Fig:UnknownPeriod_Results_M3} shows the resulting EnKF time series estimates of the seven TVP model coefficients and period $P$ when $M=3$.
Table~\ref{Tab:UnknownPeriod_Results} lists the posterior sample mean of the estimated coefficients and $P$ values for each $M$, and
Figure~\ref{Fig:UnknownPeriod_Results_TVP} shows the resulting $\bar{\theta}_{M}(t)$ approximations and time series estimates of $P$ compared with the true $\theta(t)$ and $P$, respectively, when $M=1$, 3, and 5.

Similar to the results obtained using $M=3$ with a known period (shown in Figure~\ref{Fig:KnownPeriod_Results_M3}), we see in Figure~\ref{Fig:UnknownPeriod_Results_M3} that most of the coefficient estimates converge after assimilating about one period of data; the uncertainty encoded in the $\pm2$ standard deviations around the mean is a bit wider but continues to decrease as more data are sequentially incorporated.  
The estimate of $P$ takes more time to converge but, after assimilating about two periods of data, converges closely to the true underlying period (with a relative error of approximately $8.6124\times10^{-4}$). 
While not shown, the resulting time series estimates for both position and velocity are similar to those obtained in the known period case.

The results in Table~\ref{Tab:UnknownPeriod_Results} and plots in Figure~\ref{Fig:UnknownPeriod_Results_TVP} emphasize that, again for this case, $M=3$ gives the best overall approximation to $\theta(t)$ (with smallest scaled RMSE) as well as the best estimate of $P$, with reasonably small uncertainty in this estimate by the end of the filtering process.
When $M=1$ and $M=2$, the period estimates somewhat diverge from the truth and result in under-approximations.
The results when $M=4$, while not shown graphically, are similar to $M=3$ but with more uncertainty in the posterior estimate of $P$ and an increase in the scaled RMSE of the TVP approximation.
When $M=5$, the increasing uncertainty in the posterior estimate of $P$ becomes more clear (as seen in Figure~\ref{Fig:UnknownPeriod_Results_TVP}), along with increased error in the TVP approximation.
Following from these results, the corresponding mass-spring model predictions using different initial conditions with these TVP approximation models for $\theta(t)$ are most accurate when $M=3$.
However, if directly using the models when $M=4$ or $M=5$ in this case, the increased uncertainty in the period and corresponding increased error in the TVP approximations lead to increased error in the model predictions.
This can be addressed by using the posterior mean estimate of $P$, fixing the period to this known value, and re-running the filtering process (now assuming a known period) to obtain improved TVP approximation model coefficient estimates.

\begin{table}[t!]
\caption{EnKF posterior sample means of the TVP model coefficients and period $P$ for different values of $M$ when approximating the sinusoidal forcing parameter $\theta(t) = 2\sin(t) - 0.5\cos(2t/3)$ in \eqref{Eq:MS_system}, using the data in Figure~\ref{Fig:MS_Data} and treating the period of $\theta(t)$ as an additional unknown. Numbers are rounded to four decimal places.}
\label{Tab:UnknownPeriod_Results}
\centering
\def\arraystretch{1.25}
\begin{tabular}{cccccc}
Coefficient & $M=1$ & $M=2$ & $M=3$ & $M=4$ & $M=5$ \\[0.2mm]
\hline 
$c_0$ & -0.3468  & \ 0.1705 & -0.0053  & \ 0.0962 & \ 0.0014 \\[0.2mm]
$c_1$ & \ 0.5621 & -0.0454  & -0.0043  & -0.0709  & -0.1448 \\[0.2mm]
$c_2$ & \ 0.2776 & \ 0.0085 & -0.0760  & \ 0.0274 & -0.1032 \\[0.2mm]
$c_3$ & -- & \ 1.7142 & -0.0479  & -0.0343 & \ 0.0747 \\[0.2mm]
$c_4$ & -- & -0.0627  & -0.5039  & -0.5476 & -1.3471 \\[0.2mm]
$c_5$ & -- & -- & \ 2.0026 & \ 1.9663 & \ 2.4455 \\[0.2mm]
$c_6$ & -- & -- & \ 0.0926 & -0.4897  & -0.8694  \\[0.2mm]
$c_7$ & -- & -- & -- & \ 0.0487 & \ 0.3266 \\[0.2mm]
$c_8$ & -- & -- & -- & -0.2534  & \ 0.4760 \\[0.2mm]
$c_9$ & -- & -- & -- & -- & -0.0178 \\[0.2mm]
$c_{10}$ & -- & -- & -- & -- & \ 0.3467 \\[0.2mm]
$P$ & 10.0386 & 12.5959 & 18.8658  & 18.8024 & 18.9363 \\[0.2mm]
\hline
Scaled RMSE & \ 1.0373  & \ 0.3170  & \ 0.0554  & \ 0.2288  & \ 0.8050 
\end{tabular}
\end{table}


\subsection{Example: Non-Periodic Time-Varying Parameters}

In the previous simulations, we demonstrate the effectiveness of the proposed estimation method when approximating a periodic TVP, addressing situations when we know (and fix) and when we don't know (and estimate) the underlying period.
Here we extend this approach to approximate more general TVP, in particular, considering cases when $\theta(t)$ is not periodic over the time interval of available data.
As discussed in Section~\ref{Sec:Approx_Model}, this prevents use of a direct formula for setting the $\omega_i$ values in our sinusoidal approximation model terms.
Instead, here we choose a fixed increment $\omega$ and let $\omega_i = \omega i$, $i = 1,\dots, M$, in order to systematically include sine and cosine pairs with different frequencies in the approximation.

\begin{figure}[t!] 
\centerline{\includegraphics[width = \textwidth]{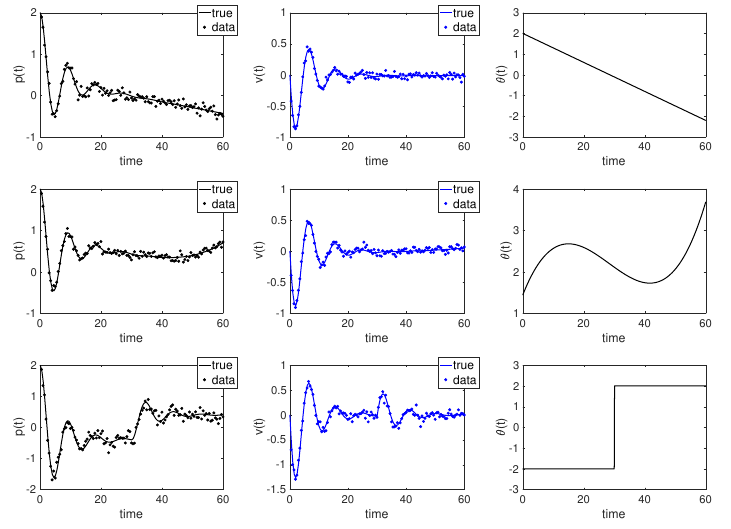}}
\caption{Position and velocity data generated from the mass-spring system in \eqref{Eq:MS_system} with the linear forcing parameter in \eqref{Eq:Linear_Forcing} (top row); cubic forcing parameter in \eqref{Eq:Cubic_Forcing} (middle row); and step forcing parameter in \eqref{Eq:Step_Forcing} (bottom row). }
\label{Fig:MS_Data_NonPeriodic}
\end{figure}

Using the same procedure as before, we simulate data from \eqref{Eq:MS_system} over the time interval [0,60] with three different non-periodic forcing parameters: 
\begin{itemize}
\item[(i)] a linear polynomial, where 
\begin{equation}\label{Eq:Linear_Forcing}
\theta(t) = -0.07t + 2;
\end{equation}
\item[(ii)] a cubic polynomial, where
\begin{equation}\label{Eq:Cubic_Forcing}
\theta(t) = 0.0001(t - 25)^3 - 0.001t^2 + 3;
\end{equation}
and 
\item[(iii)] a step function, where
\begin{equation}\label{Eq:Step_Forcing}
\theta(t) = \begin{cases} -2 & t \leq 30 \\ \ \ 2 & t>30  \end{cases} .
\end{equation}
\end{itemize}
Figure~\ref{Fig:MS_Data_NonPeriodic} shows the simulated data in each case.

\begin{figure}[t!] 
\centerline{\includegraphics[width = \textwidth]{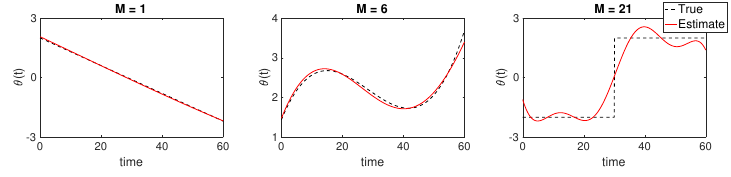}}
\caption{Resulting Fourier series-based approximations of three non-periodic forcing parameters $\theta(t)$ in \eqref{Eq:MS_system} for the values of $M$ yielding the smallest scaled RMSE, using the data in Figure~\ref{Fig:MS_Data_NonPeriodic} and setting $\omega_i = 0.01 i$ for $i=1,\dots,M$ in the TVP approximation models.
From left to right: the linear polynomial in \eqref{Eq:Linear_Forcing}; the cubic polynomial in \eqref{Eq:Cubic_Forcing}; and the step function in \eqref{Eq:Step_Forcing}.
On each plot, the solid red line denotes the TVP approximation $\bar{\theta}_{M}(t)$ computed using the posterior mean coefficient values and the dashed black line shows the true underlying $\theta(t)$.
}
\label{Fig:NonPeriodic_Results_TVP}
\end{figure}

In building the TVP approximation models as in \eqref{Eq:TV_form2}, we set an increment of $\omega = 0.01$ and define $\omega_i = 0.01 i$ for each $i = 1,\dots, M$.
Since each of the underlying $\theta(t)$ have different levels of complexity, we test a variety of $M$ values in each case.
Figure~\ref{Fig:NonPeriodic_Results_TVP} shows the best resulting approximations (i.e., those with the smallest corresponding scaled RMSEs) for each non-periodic forcing parameter.
For the linear forcing parameter, the best fit occurs when $M=1$ (with scaled RMSE $\approx$ 0.0239), needing only three terms in the TVP approximation model to obtain an accurate approximation.  
Here, increasing $M$ to somewhat larger values (e.g., $M=4$ or 5) increases the approximation error and number of unknowns but still permits reasonable approximations with convergent coefficients.

More terms are needed when approximating the cubic forcing parameter, and it is important to note that, for this example, using too small an $M$ value results in coefficient estimates that do not converge over the time interval of available data.
More specifically, the coefficients do not converge when $M=1$, 2 or 3, but convergence improves beginning with $M=4$.
The best fit for the cubic forcing parameter occurs when $M=6$ (with scaled RMSE $\approx$ 0.1596), but with a similarly good fit when $M=5$ (scaled RMSE $\approx$ 0.1614) requiring the estimation of two less coefficients.

Estimating the step forcing parameter presents the most difficult challenge of the three, given the jump discontinuity in the function halfway through the time interval of available data.
For this example, the TVP approximation model coefficients do not converge for smaller $M$ values, but a significant increase in $M$ leads to reasonable fits with convergent coefficients. 
The TVP approximations yield similar scaled RMSE values for $M$ between 15 and 22, with the lowest occurring when $M=21$ (scaled RMSE $\approx$ 0.2737).
While not capturing the exact shape, the TVP approximations in this case are able to capture the jump point between constant parameter values in the underlying step function.
The mass-spring model predictions for different initial conditions follow as expected using the best TVP approximation results over the time interval [0,60], but it becomes more difficult to accurately predict the behavior of the system after this time frame for non-periodic TVPs, since the parameters themselves will continue to dynamically change over time without updating their approximation models.

\section{Conclusions and Future Work}
\label{Sec:Conclusions}

In this work, we present a novel Fourier series-based approximation method for estimating time-varying parameters in deterministic dynamical systems, where we represent $\theta(t)$ as a linear combination of sine and cosine functions and estimate the unknown approximation model coefficients using ensemble Kalman filtering.
This approach allows us to construct accurate TVP approximations $\bar{\theta}_{M}(t)$, dependent on the integer parameter $M$ (which dictates the number of approximation model terms) and posterior EnKF mean estimates for the model coefficients.
With several numerical examples using a forced mass-spring system, we illustrate the effectiveness of this approach in approximating TVPs of different forms, including cases when $\theta(t)$ is periodic with a known period, periodic with an unknown period (which is jointly estimated), and non-periodic over the time interval of observed data.

One important aspect to consider in successfully applying this method is how to best select $M$ in advance of running the estimation procedure.
Our goal in practice is to choose the smallest integer $M$ that will result in a reasonable TVP approximation and corresponding ODE system predictions, in order to keep the number of unknowns to a minimum as well as to reduce approximation errors from including additional terms that may not be needed.
Another important consideration is how to appropriately set the angular frequencies $\omega_i$, $i=1,\dots,M$, of the sinusoidal terms in the TVP approximation models.
When $\theta(t)$ is periodic, we are able to set $\omega_i$ directly using the formula in \eqref{Eq:Known_Omega} and either fix or estimate $P$, depending on what information is available.  
However, when estimating TVPs that are not periodic (or not known to be periodic), we select a fixed increment $\omega$ and set $\omega_i = \omega i$; the choice of this increment is vital in the resulting TVP approximation and also can affect which $M$ is most appropriate for the problem.
Future work will include incorporating model selection techniques \cite{Stoica2004, Ding2018} and potential pre-processing steps to select $M$ and $\omega$ systematically for a given problem.

We note that the proposed approach is not limited to the use of EnKF in estimating the unknown coefficients and could also be implemented using different nonlinear filtering methods (e.g., particle filters \cite{Kantas2015}) with sequentially-arriving data or non-sequential Bayesian approaches (e.g., MCMC methods \cite{Luengo2020}) if the full time series of data is available at once.  
However, our results show that we are able to reasonably approximate the unknown TVPs of interest using the EnKF with a relatively small sample size compared to those generally needed for particle filtering or MCMC-based approaches.  
While the examples in this work each focus on estimating a single, univariate TVP for the system considered, future work will examine the feasibility of this approach in simultaneously estimating multiple TVPs for a given system, as well as introducing model-data mismatch as we move toward real-data application.



\section*{Acknowledgements}

This work was supported by the National Science Foundation under grant number NSF/DMS-1819203 (A. Arnold).


\section*{ORCID iDs}
Andrea Arnold: \url{https://orcid.org/0000-0003-3003-882X}


\begin{appendices}

\renewcommand\thefigure{\thesection.\arabic{figure}} 
\setcounter{figure}{0}

\renewcommand\thetable{\thesection.\arabic{table}} 
\setcounter{table}{0}

\renewcommand\theequation{\thesection.\arabic{equation}} 
\setcounter{equation}{0}


\section{Numerical Results: Example Using Derivative-Based Augmentation}
\label{Appendix:Augmentation_Approach}

\begin{figure}[t!] 
\centerline{\includegraphics[width = \textwidth]{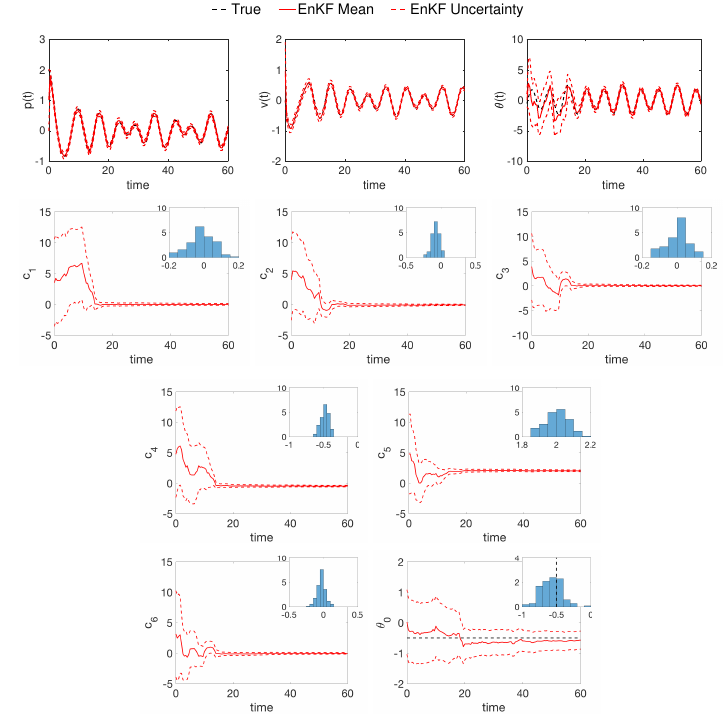}}
\caption{EnKF time series estimates of position, velocity, and the forcing parameter $\theta(t)$, along with the six unknown coefficients in \eqref{Eq:TV_deriv} and the initial condition $\theta_0$ of the sinusoidal forcing parameter $\theta(t) = 2\sin(t) - 0.5\cos(2t/3)$ in \eqref{Eq:MS_system} using the derivative-based augmentation approach when $M=3$, given the data in Figure~\ref{Fig:MS_Data} and assuming a known period for $\theta(t)$. 
On each plot, the solid red line denotes the EnKF sample mean and the dashed red lines show $\pm2$ standard deviations around the mean.
The plots for each coefficient, $c_1, \dots, c_6$, and for $\theta_0$ also show the resulting histograms of the posterior samples at time $t=60$. 
}
\label{Fig:Augmented_KnownPeriod_Results_M3}
\end{figure}

\begin{figure}[ht!] 
\centerline{\includegraphics[width = \textwidth]{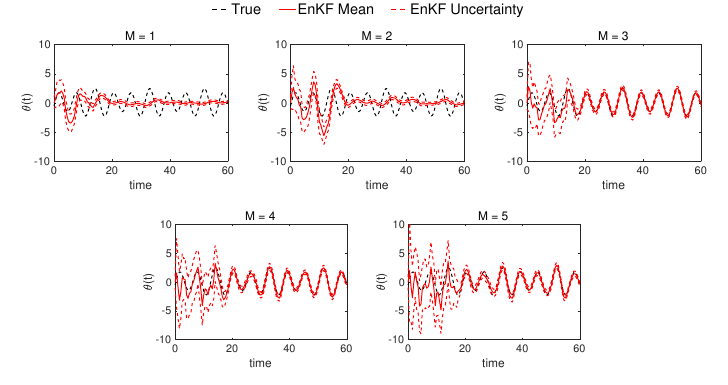}}
\caption{EnKF time series estimates of the sinusoidal forcing parameter $\theta(t) = 2\sin(t) - 0.5\cos(2t/3)$ in \eqref{Eq:MS_system} when using the derivative-based augmentation approach with different values of $M$, given the data in Figure~\ref{Fig:MS_Data} and assuming a known period for $\theta(t)$.
On each plot, the solid red line denotes the EnKF sample mean, the dashed red lines show $\pm2$ standard deviations around the mean, and the dashed black line shows the true underlying $\theta(t)$.
}
\label{Fig:Augmented_KnownPeriod_Results_TVP}
\end{figure}

In this section, we apply the derivative-based augmentation approach described in Section~\ref{Sec:Alternative} to estimate the sinusoidal forcing parameter $\theta(t) = 2\sin(t) - 0.5\cos(2t/3)$ in \eqref{Eq:MS_system} using the same data as in Figure~\ref{Fig:MS_Data} and assuming a known period for $\theta(t)$.
Compared to the approach used in Section~\ref{Sec:KnownPeriod_Ex}, the difference with this implementation is that we now treat $\theta(t)$ as unobserved system state, where
\begin{equation}
\frac{d\theta}{dt} = \theta'_M(t), \quad \theta(0) = \theta_0
\end{equation}
with $\theta'_M(t)$ defined as in \eqref{Eq:TV_deriv}, and we track it along with $p(t)$ and $v(t)$, estimating the $2M$ coefficients $c_1,\dots,c_{2M}$ plus the TVP initial value $\theta_0$ for a total of $2M+1$ unknown constant parameters.
Figure~\ref{Fig:Augmented_KnownPeriod_Results_M3} shows the resulting EnKF time series estimates of position, velocity, and $\theta(t)$, along with the six unknown coefficients in \eqref{Eq:TV_deriv} and the initial condition $\theta_0$, when $M=3$.
Table~\ref{Tab:Augmented_KnownPeriod_Results} lists the posterior sample mean of each coefficient and $\theta_0$ for each $M$ value, and
Figure~\ref{Fig:Augmented_KnownPeriod_Results_TVP} shows the corresponding EnKF time series estimates of $\theta(t)$ in each case compared with the true $\theta(t)$.

The results in Figure~\ref{Fig:Augmented_KnownPeriod_Results_M3} highlight the differences between this derivative-based augmented systems implementation and the direct approximation model (with corresponding results shown in Figure~\ref{Fig:KnownPeriod_Results_M3}), where here $\theta(t)$ is treated as an unobserved system state and tracked along with position and velocity.
The six coefficients converge after assimilating about one period of data, which is also reflected in the tracking of $\theta(t)$, while the estimate of $\theta_0$ continues to improve over the remaining time interval of observations. 
As the plots in Figure~\ref{Fig:Augmented_KnownPeriod_Results_TVP} illustrate, the filter is not able to well track $\theta(t)$ when $M=1$ and 2, noticeably under-estimating the dynamics with a damping effect beginning close to the one period mark (i.e., around time $t=19$).
The estimation significantly improves when $M=3$, where the EnKF mean estimate very well captures the behavior of the true underlying parameter after assimilating one period of data, with fairly tight uncertainty bounds around the mean estimate. 
Similar behavior occurs when $M=4$ and 5, although there is more initial error in the approximation over the first period and slightly wider uncertainty bounds after the coefficients converge.
The results in Table~\ref{Tab:Augmented_KnownPeriod_Results} paint a similar picture for the posterior mean estimates of $\theta_0$, where the estimate when $M=3$ yields the smallest relative error ($\approx$ 0.1440) when compared to the true initial value $\theta_0 = -0.5$.

\begin{table}[ht!]
\caption{EnKF posterior sample means of the TVP model coefficients and initial value $\theta_0$ for different values of $M$ when using the derivative-based augmentation approach to estimate the sinusoidal forcing parameter $\theta(t) = 2\sin(t) - 0.5\cos(2t/3)$ in \eqref{Eq:MS_system} given the data in Figure~\ref{Fig:MS_Data} and assuming a known period for $\theta(t)$. Numbers are rounded to four decimal places.}
\label{Tab:Augmented_KnownPeriod_Results}
\centering
\def\arraystretch{1.25}
\begin{tabular}{cccccc}
Coefficient & $M=1$ & $M=2$ & $M=3$ & $M=4$ & $M=5$ \\[0.2mm]
\hline 
$c_1$ & \ 0.1410  & -0.0439  & -0.0134 & \ 0.0590  & -0.0324  \\[0.2mm]
$c_2$ & -0.0073  & \ 0.1723 & -0.0814 & -0.0612 & -0.0894  \\[0.2mm]
$c_3$ & -- & -0.0865 & \ 0.0034  & -0.0505 & \ 0.0073  \\[0.2mm]
$c_4$ & -- & -0.1852 & -0.4840  & -0.4655 & -0.5573 \\[0.2mm]
$c_5$ & -- & -- & \ 2.0026 & \ 2.0152  & \ 1.9659 \\[0.2mm]
$c_6$ & -- & -- & -0.0318 & \ 0.0115  & \ 0.0038 \\[0.2mm]
$c_7$ & -- & -- & -- & \ 0.1162 & \ 0.1577 \\[0.2mm]
$c_8$ & -- & -- & -- & \ 0.1301 & \ 0.1571 \\[0.2mm]
$c_9$ & -- & -- & -- & -- & -0.2124 \\[0.2mm]
$c_{10}$ & -- & -- & -- & -- & \ 0.1677  \\[0.2mm]
$\theta_0$ & -0.0643  & \ 0.3544  & -0.5720  & -0.4168  & -0.2902
\end{tabular}
\end{table}

\end{appendices}



\bibliography{paper_refs}{}

\begin{thebibliography}{10}

\bibitem{Sika2008}
G.~Sika and P.~Velex.
\newblock Instability analysis in oscillators with velocity-modulated
  time-varying stiffness--applications to gears submitted to engine speed
  fluctuations.
\newblock {\em Journal of Sound and Vibration}, 318:166--175, 2008.

\bibitem{Shen2006}
Y.~Shen, S.~Yang, and X.~Liu.
\newblock Nonlinear dynamics of a spur gear pair with time-varying stiffness
  and backlash based on incremental harmonic balance method.
\newblock {\em International Journal of Mechanical Sciences}, 48:1256--1263,
  2006.

\bibitem{Batouli2015}
J.~Batouli, M.~El~Baz, and A.~Maaouni.
\newblock {RLC} circuit realization of a q-deformed harmonic oscillator with
  time dependent mass.
\newblock {\em Physics Letters A}, 379:1619--1626, 2015.

\bibitem{Preiner2007}
J.~Preiner, J.~Tang, V.~Pastushenko, and P.~Hinterdorfer.
\newblock Higher harmonic atomic force microscopy: Imaging of biological
  membranes in liquid.
\newblock {\em Physical Review Letters}, 99:046102, 2007.

\bibitem{Altizer2006}
S.~Altizer, A.~Dobson, P.~Hosseini, P.~Hudson, M.~Pascual, and P.~Rohani.
\newblock Seasonality and the dynamics of infectious diseases.
\newblock {\em Ecology Letters}, 9(467-484), 2006.

\bibitem{Zeng2020}
Y.~Zeng, X.~Guo, Q.~Deng, S.~Luo, and H.~Zhang.
\newblock Forecasting of {COVID-19}: Spread with dynamic transmission rate.
\newblock {\em Journal of Safety Science and Resilience}, 1(2):91--96, 2020.

\bibitem{Calvetti2020}
D.~Calvetti, A.~P. Hoover, J.~Rose, and E.~Somersalo.
\newblock Metapopulation network models for understanding, predicting, and
  managing the coronavirus disease {COVID-19}.
\newblock {\em Frontiers in Physics}, 8:261, 2020.

\bibitem{Linaro2018}
D.~Linaro, I.~Biro, and M.~Giugliano.
\newblock Dynamical response properties of neocortical neurons to
  conductance-driven time-varying inputs.
\newblock {\em European Journal of Neuroscience}, 47:17--32, 2018.

\bibitem{Shamir2009}
M.~Shamir, O.~Ghitza, S.~Epstein, and N.~Kopell.
\newblock Representation of time-varying stimuli by a network exhibiting
  oscillations on a faster time scale.
\newblock {\em PLOS Computational Biology}, 5(5):e1000370, 2009.

\bibitem{Campbell2020}
K.~Campbell, L.~Staugler, and A.~Arnold.
\newblock Estimating time-varying applied current in the {Hodgkin-Huxley}
  model.
\newblock {\em Applied Sciences}, 10(2):550, 2020.

\bibitem{ArnoldFichera2022}
A.~Arnold and L.~Fichera.
\newblock Identification of tissue optical properties during thermal
  laser-tissue interactions: an ensemble {Kalman} filter-based approach.
\newblock {\em International Journal for Numerical Methods in Biomedical
  Engineering}, 38(4):e3574, 2022.

\bibitem{Bashkatov2016}
A.~N. Bashkatov, E.~A. Genina, V.~I. Kochubey, and V.~V. Tuchin.
\newblock Quantification of tissue optical properties: perspectives for precise
  optical diagnostics, phototherapy and laser surgery.
\newblock {\em Journal of Physics D: Applied Physics}, 49(50):501001, 2016.

\bibitem{Labiod2013}
S.~Labiod, H.~Boubertakh, and T.~M. Guerra.
\newblock Fourier series-based adaptive tracking control for robot
  manipulators.
\newblock In {\em Proceedings of the 3rd International Conference on Systems
  and Control}, pages 968--972. IEEE, 2013.

\bibitem{Caruso2021}
A.~Caruso, M.~Bassetto, G.~Mengali, and A.~A. Quarta.
\newblock Optimal solar sail trajectory approximation with finite {Fourier}
  series.
\newblock {\em Advances in Space Research}, 67:2834--2843, 2021.

\bibitem{Konidaris2011}
G.~Konidaris, S.~Osentoski, and P.~Thomas.
\newblock Value function approximation in reinforcement learning using the
  {Fourier} basis.
\newblock In {\em Proceedings of the Twenty-Fifth AAAI Conference on Artificial
  Intelligence}, pages 380--385, 2011.

\bibitem{Liuzzo2007}
S.~Liuzzo, R.~Marino, and P.~Tomei.
\newblock Adaptive learning control of nonlinear systems by output error
  feedback.
\newblock {\em IEEE Transactions on Automatic Control}, 52(7):1232--1248, 2007.

\bibitem{Chen2010}
W.~Chen, W.~Li, and Q.~Miao.
\newblock {Backstepping control for periodically time-varying systems using
  high-order neural network and Fourier series expansion}.
\newblock {\em ISA Transactions}, 49:283--292, 2010.

\bibitem{Zhang2012}
C.-L. Zhang and J.-M. Li.
\newblock {Hybrid function projective synchronization of chaotic systems with
  uncertain time-varying parameters via Fourier series expansion}.
\newblock {\em International Journal of Automation and Computing}, 9:388--394,
  2012.

\bibitem{Chen2021}
J.~Chen and J.~Li.
\newblock Distributed consensus control of periodically time-varying
  multi-agent systems using neural networks and fourier series expansion.
\newblock {\em Journal of the Franklin Institute}, 358:7170--7186, 2021.

\bibitem{Arnold2018}
A.~Arnold and A.~L. Lloyd.
\newblock An approach to periodic, time-varying parameter estimation using
  nonlinear filtering.
\newblock {\em Inverse Problems}, 34:105005, 2018.

\bibitem{Arnold2020}
A.~Arnold.
\newblock Using {Monte Carlo} particle methods to estimate and quantify
  uncertainty in periodic parameters.
\newblock In {\em Advances in Mathematical Sciences}, pages 213--226. Springer,
  2020.

\bibitem{Feng2017}
M.~Feng, P.~Liu, S.~Guo, L.~Shi, C.~Deng, and Ming B.
\newblock Deriving adaptive operating rules of hydropower reservoirs using
  time-varying parameters generated by the {EnKF}.
\newblock {\em Water Resources Research}, 53:6885--6907, 2017.

\bibitem{Voss2004}
H.~U. Voss, J.~Timmer, and J.~Kurths.
\newblock Nonlinear dynamical system identification from uncertain and indirect
  measurements.
\newblock {\em International Journal of Bifurcation and Chaos}, 14:1905--1933,
  2004.

\bibitem{Bian2011}
X.~Bian, X.~R. Li, H.~Chen, D.~Gan, and J.~Qiu.
\newblock Joint estimation of state and parameter with synchrophasors -- {Part
  II}: parameter tracking.
\newblock {\em IEEE Transactions on Power Systems}, 26:1209--1220, 2011.

\bibitem{Arnold2019}
A.~Arnold.
\newblock Exploring the effects of uncertainty in parameter tracking estimates
  for the time-varying external voltage parameter in the {FitzHugh-Nagumo}
  model.
\newblock In {\em 6th International Conference on Computational and
  Mathematical Biomedical Engineering}, pages 512--515, 2019.

\bibitem{Calvetti2021}
D.~Calvetti, A.~Hoover, J.~Rose, and E.~Somersalo.
\newblock Bayesian particle filter algorithm for learning epidemic dynamics.
\newblock {\em Inverse Problems}, 37:115008, 2021.

\bibitem{Arnold2023}
A.~Arnold.
\newblock When artificial parameter evolution gets real: particle filtering for
  time-varying parameter estimation in deterministic dynamical systems.
\newblock {\em Inverse Problems}, 39:014002, 2023.

\bibitem{Evensen1994}
G.~Evensen.
\newblock {Sequential data assimilation with a nonlinear quasi‐geostrophic
  model using Monte Carlo methods to forecast error statistics}.
\newblock {\em Ocean Dynamics}, 99:10143--10162, 1994.

\bibitem{Burgers1998}
G.~Burgers, P.~J. van Leeuwen, and G.~Evensen.
\newblock Analysis scheme in the ensemble {Kalman} filter.
\newblock {\em Monthly Weather Review}, 126:1719--1724, 1998.

\bibitem{Arnold2014}
A.~Arnold, D.~Calvetti, and E.~Somersalo.
\newblock {Parameter estimation for stiff deterministic dynamical systems via
  ensemble Kalman filter}.
\newblock {\em Inverse Problems}, 30:105008, 2014.

\bibitem{Katzfuss2016}
M.~Katzfuss, J.~R. Stroud, and C.~K. Wikle.
\newblock Understanding the ensemble {Kalman} filter.
\newblock {\em The American Statistician}, 70:350--357, 2016.

\bibitem{Evensen2009}
G.~Evensen.
\newblock The ensemble {Kalman} filter for combined state and parameter
  estimation.
\newblock {\em IEEE Control Systems Magazine}, 29:83--104, 2009.

\bibitem{Shampine2008}
L.~F. Shampine.
\newblock Vectorized adaptive quadrature in {MATLAB}.
\newblock {\em Journal of Computational and Applied Mathematics}, 211:131--140,
  2008.

\bibitem{Nagle2011}
R.~K. Nagle, E.~B. Saff, and A.~D. Snider.
\newblock {\em Fundamentals of Differential Equations and Boundary Value
  Problems}.
\newblock Pearson, 6 edition, 2011.

\bibitem{Boyce2001}
W.~E. Boyce and R.~C. DiPrima.
\newblock {\em Elementary Differential Equations and Boundary Value Problems}.
\newblock John Wiley \& Sons, New York, 7th edition, 2001.

\bibitem{Stoica2004}
P.~Stoica and Y.~Selen.
\newblock Model-order selection: a review of information criterion rules.
\newblock {\em IEEE Signal Processing Magazine}, 21(4):36--47, 2004.

\bibitem{Ding2018}
J.~Ding, V.~Tarokh, and Y.~Yang.
\newblock Model selection techniques: An overview.
\newblock {\em IEEE Signal Processing Magazine}, 35(6):16--34, 2018.

\bibitem{Kantas2015}
N.~Kantas, A.~Doucet, S.~S. Singh, J.~Maciejowski, and N.~Chopin.
\newblock On particle methods for parameter estimation in state-space models.
\newblock {\em Statistical Science}, 30(3):328--351, 2015.

\bibitem{Luengo2020}
D.~Luengo, L.~Martino, M.~Bugallo, V.~Elvira, and S.~Sarkka.
\newblock A survey of {Monte Carlo} methods for parameter estimation.
\newblock {\em EURASIP Journal on Advances in Signal Processing}, 2020:25,
  2020.

\end{thebibliography}

\end{document}